\documentclass[12pt]{amsart}

\usepackage{amsmath, amssymb}
\usepackage[margin=7.5em]{geometry}
\usepackage{enumitem}

\newcommand{\A}{\mathbb{A}1}
\newcommand{\DD}{\mathcal{D}}
\newcommand{\DDD}{{\mathcal{D}_{\A1}}}
\newcommand{\dd}{\partial}

\newcommand{\Q}{\mathbb{Q}}
\newcommand{\Z}{\mathbb{Z}}

\newcommand{\C}{\mathbb{C}}

\renewcommand{\P}{\mathbb{P}}

\newcommand{\R}{\mathbb{R}}

\newcommand{\Frob}{\text{Frob}}

\newcommand{\Gal}{\mathrm{Gal}}
\DeclareMathOperator{\convo}{{\ast \,}}

\def\Gm{{\mathbf{G}}_m} 
\def\DD{{\mathcal{D}}}

\newcommand{\mathscr}{\mathcal}

\newcommand{\mystrut}{{\vrule depth 1em width 0em height 2em}}

\def\A1{\mathbb{A}^1}
\def\P1{{{\Bbb{P}}^1}}

\def\Gm{{\bf G_m}}

\def\SL2{{\mathrm SL2}}

\renewcommand{\P}{{\Bbb{P}}}

\newcounter{pphcounter}[section]
 \renewcommand{\thepphcounter}{\arabic{pphcounter}}
\newcommand{\pph}{\medskip \refstepcounter{pphcounter}
 \bf \thepphcounter.~\rm}

\newcounter{proofstepcounter}[pphcounter]
\renewcommand{\theproofstepcounter}{\arabic{proofstepcounter}}
\newcommand{\proofstep}{\medskip \refstepcounter{proofstepcounter}
	\hspace{0.125em} \sl \theproofstepcounter.~\rm}

\def \DD {\mathcal{D}}

\def\OO{{\mathcal{O}}}
\renewcommand{\phi}{{\varphi}}

\def\Alpha{{\mathrm{{\mathrm{A}}}}}
\def\Beta{{\mathrm{B}}}

\newcommand\Exp[1]{e^{#1}}

\newcommand{\FT}{{\mathrm{FT}}}
\newcommand{\VM}{{\mathrm{V}}}

\def\isomorphic{\cong}

\newcommand{\tensor}{\otimes}

\newcommand{\ta}{{\tilde \alpha}}

\newcommand{\tb}{{\tilde \beta}}

\renewcommand{\FT}{\mathop{\mathrm{FT}}}
\newcommand{\ii}{\mathsf{i}}
                                           
\begin{document}

\title{On $L$--derivatives and  biextensions \\ 
of Calabi--Yau motives}
\author{V. Golyshev}

\begin{abstract}
\noindent
We prove that certain differential operators of the form $DLD$  with $L$  hypergeometric and $D=z \frac{d}{dz}$ are of Picard--Fuchs type.  We  give closed hypergeometric expressions for minors of the biextension period matrices that arise from certain rank 4 weight 3 Calabi--Yau motives  presumed to be of analytic rank 1. We compare their values numerically to the first derivative of the $L$--functions of the respective motives at 
$s=2$. 
\end{abstract}
\maketitle

\bigskip

The goal of this note is to explain, favoring expedience over detail, how one can systematically obtain  explicit numerical evidence in support of a Birch--Swinnerton-Dyer-type conjecture for hypergeometric Calabi--Yau motives.
For a Calabi--Yau threefold $X/\Q$ with
Hodge numbers $h^{3,0} = h^{0,3} =1,\: h^{2,1} = h^ {1,2} = a$,
Poincar\'e duality defines a non-degenerate
alternating form on the third cohomology $H^3(X)$ for any Weil cohomology theory.  We view the collection of such cohomological realizations as arising
from a  so--called {\em symplectic motive} of rank $2+2a$.
We will focus on the case $a=1$ here; these motives, now colloquially called `(1,1,1,1)--motives', are expected to exist in 1--parameter families \cite{GolyshevvanStraten2023}. Their typical  Euler factors can be obtained as 
characteristic polynomials of the geometric $p$--Frobenius acting in the $l$--adic cohomology of $X$ over the algebraic closure. They take
the form
\[\det(1-T \cdot \Frob_p \mid_{H^3_\text{\'et}(\bar X,\Q_l)})=1+\alpha_pT+\beta_p p T^2+ p^3\alpha_p T^3+p^6T^4\]
with $\alpha_p, \beta_p \in \Z$.

It is believed that the completed 
$L$--function 
\[\Lambda(s)=\left(\frac{N}{\pi^4}\right)^{s/2} \Gamma\left(\frac{s-1}{2}\right)\Gamma\left(\frac{s}{2}\right)\Gamma\left(\frac{s}{2}\right)\Gamma\left(\frac{s+1}{2}\right) L(H^3(X),s), \]
is entire and satisfies $\Lambda(s) = \pm \, \Lambda(4-s)$, where $N$ is the conductor. The known meromorphicity and existence of a functional equation \cite{PatrikisTaylor2015} enable one in principle to study the leading coefficient of the Taylor series of $L(s)$ along the lines suggested by the conjectures of Deligne \cite{Deligne1979} and Birch--Swinnerton-Dyer (and Beilinson \cite{Beilinson1987}, Bloch \cite{Bloch1980, Bloch1984}, Gillet--Soul\'e \cite{GilletSoule1984}). More broadly, if one is to think of these $(1,1,1,1)$--motives as analogues of elliptic curves over $\Q$ two dimensions higher, a question arises of what standard motivic conjectures known to be true or confirmed numerically for elliptic curves survive in this new setup. The key and probably indispensable ingredient here will be a suitable automorphy theorem. Many believe, for instance, that a weight 3 paramodular newform (a Hecke--eigen (3,0)--regular form on the Siegel threefold parametrizing $(1,N)$--polarized abelian surfaces) $f_M$ could be assocated to such a motive $M$ of conductor $N$ so that $L(f_M,s)= L(M,s)$. 
With automorphy proven --- in general, or for any given 
motive~$M$ --- one could try to proceed by relating the central $L$--value (or the leading coefficient) at $s=2$ obtained from an {integral representation} for the $L$--function to a certain Hodge--theoretic volume arising in a biextension of $M$, an idea that can be traced back to Bloch's early work \cite{Bloch1980}; see also  \cite{BlochKato1990, KontsevichZagier2001, Scholl1991}.

In analytic rank $1$, one would seek a $\mathrm{GSp}(4)$--analogue of the Gross--Zagier formula \cite{GrossZagier1986} that might express $L'(2,M)$ in terms of the height pairing between certain curves on the Siegel threefold parametrizing special abelian surfaces.  Its proof, however, is expected to be very difficult and not to be found soon, so a numerical study is desirable as a second-best. Once the Dirichlet series of $M$ and the shape of the functional equation are known, the technology described in \cite{Dokchitser2004}
and implemented in Magma \cite{BosmaCannonPlayoust1993}, enables  one (in principle) to compute the Taylor expansion of $L (M,s)$ to an arbitrary precision. 

The paper \cite{RobertsVillegas2021} is an excellent introduction to hypergeometric motives and explains how to compute hypergeometric $L$--functions. The present note can be viewed by the reader as a companion paper. We show how a combination of two ideas specific to hypergeometric pencils enables one to write down closed formulas for the (archimedean) extension volumes and obtain evidence in support of B-SD. One is the principle that gamma structures \cite{GolyshevMellit2014} give rise to Betti structures. The other says that the motive of the \emph{total space} 
of a hypergeometric pencil can be used to provide every fiber with a biextension Hodge structure \cite{BlochdeJongSertoz2022}.
The relevant biextension can be viewed as joining together two Katz's extensions \cite[8.4.7, 8.4.9]{Katz1990} going the opposite directions.

\pph \bf Hypergeometric (1,1,1,1)--families. \rm  
The arithmetic of some of these hypergeometric familes was studied by e.g. Dwork \cite{Dwork1969} and Schoen \cite{Schoen1986}. 
The interest in families of $\Q$--Calabi--Yau motives with points of maximally  
unipotent monodromy surged in the wake of the discovery of mirror symmetry 
\cite{CandelasOssaGreenParkes1991}. The simplest are the $14$ hypergeometric families, which directly
generalize the famous Dwork pencil \cite{HofmannvanStraten2015}. These (and certain `quadratic twists' of these, as we will see) are probably the most amenable to direct computation with the $l$-adic and Betti-de Rham realizations.

N. Katz introduced implicitly the concept of a hypergeometric motivic sheaf in \cite{Katz1990} by analyzing in detail hypergeometric differential equations, i.e. scalar differential equations of the form $$L_{\alpha, \beta} S(z)=0   \qquad \eqno{(*)} $$  with
$$ 
L_{\alpha, \beta}=  \prod_{i=1}^n (D-\alpha_i)- \lambda
z\prod_{j=1}^n(D-\beta_j), \;\;\; D= z\frac{d}{dz},$$ 
and proving a theorem that states that an irreducible regular singular hypergeometric
differential equation with rational indices (and $\lambda \in \bar \Q$) is motivic, i.e. arises in a piece of relative cohomology in a pencil of algebraic varieties defined over a number field. An analogue of Katz's theorem holds for tame hypergeometric $l$--adic sheaves over $\Gm / \,{\overline{\Q}}$ whose local inertiae act quasiunipotently. \nocite{Dutka1984}
If one furthermore requires that the sets
$\exp(2\pi \ii \alpha_i)$'s and $\exp(2\pi \ii \beta_j)$'s are each $\Gal(\overline{\Q}/\Q)$--stable and $\lambda \in \Q$,  a motivic construction can be defined over $\Q$, cf. \cite{BeukersCohenMellit2015}.

\pph {\bf Gamma structures give rise to Betti structures.}\label{assumptions-1}
In order to refine hypergeometric $D$--modules to Hodge modules one needs to identify the $\Q$--bases  of the spaces of local solutions that represent the periods of relative $\Q$--de Rham forms along $\Q$--Betti cycles.  
Following Dwork, one can think of hypergeometric families as deformations of Fermat hypersurfaces (with their relatively simple motivic structures) obtained by introducing an extra monomial to the defining equation. From this perspective it is clear that the leading expansion coefficients of $\Q$--Betti solutions  of hypergeometric Hodge modules should be proportional to products of the values of the gamma function at rational arguments corresponding to the hypergeometric indices.  A theorem on hypergeometric monodromy in \cite{GolyshevMellit2014}  says, in particular, the following.
Assume that (\bf A1\rm):
\begin{itemize}
\item 
the sets
$\exp(2\pi \ii \alpha_i)$'s and $\exp(2\pi \ii \beta_j)$'s are each $\Gal(\overline{\Q}/\Q)$--stable and $\lambda \in \Q$	
	
\item
$\alpha_i \neq \beta_{i'} \mod \Z$ for all $i, i'$ 

and, merely to make our statement simpler, that 

\item $\alpha_i \neq \alpha_{i'} \mod \Z$ for all $i \neq i'$ either.
\end{itemize}
 
Put
$$
\mathbf{\Gamma}(s)\,=\, \mathbf{\Gamma}_{\alpha, \beta}(s) \,=\, {\mystrut \prod_{i=1}^n \Gamma(s  - \alpha_i + 1)^{-1} \prod_{i=1}^n \Gamma (-s + \beta_i + 1)^{-1}} \quad (s\in\C),
$$
and $\Alpha_i =\Exp{2\pi\ii \alpha_i}, \;\Beta_j =\Exp{2\pi\ii \beta_j}$ To simplify notation, assume until the end of this paragraph that $\lambda = (-1)^n$. In general, $(*)$ comes with a gamma structure that
is defined to be the set $\boldsymbol\gamma = \{ \sum_{s \in s_0 +\Z} \mathbf{\Gamma}(s)\, z^s  \mid s_0 \in \C \}$ of formal solutions  to $(*)$ and is meant to specialize to a Betti structure when the hypergeometric indices are rational. In particular, consider the basis of local solutions of $(*)$ at $0$ given by
$$
S_{\Alpha_j}(z) = \sum_{l=0}^\infty  \mathbf\Gamma(l + \alpha_j)\, z^{l + \alpha_j} \, \in \, \boldsymbol\gamma .
$$
Then the monodromy of $(*)$ around $0$ is given by 
 $$M_0 (S_{\Alpha_1}(z), \dots, S_{\Alpha_n} (z))^t = (\Alpha_1S_{\Alpha_1}(z), \dots, \Alpha_n S_{\Alpha_n}(z))^t.$$  
Denote by $\VM_{\mathrm{A}}$ the respective Vandermonde matrix 
$$
\VM_{\mathrm{A}} = \begin{pmatrix}
1 & {\mathrm{A}}_1 & \cdots & {\mathrm{A}}_1^{n-1} \\
1 & {\mathrm{A}}_2  & \cdots &  {\mathrm{A}}_2^{n-1} \\
\vdots & \vdots & & \vdots
\end{pmatrix} 
.$$
The 
global monodromy of $\, (*)\,$ in the basis
$ 
V_{\mathrm{A}}^t(S_{{\mathrm{A}}_{1}}( z), \dots, S_{\mathrm{A}_n} ( z))^t$ is shown in \cite{GolyshevMellit2014} 
to be in $\mathrm{GL}_n(\Q)$, and in fact defines a $\Q$--local system that underlies a Hodge module.  

\pph \label{assumptions-2} To identify the Hodge filtration, we proceed as follows.
For simplicity, let us further assume, as is the case with our hypergeometric $(1,1,1,1)$--motives, that (\bf A2\rm):
\begin{itemize}
\item $n$ is divisible by $4$;
\item the sets $\Alpha$'s and $\Beta$'s are 
\emph{maximally non-interlaced} on the unit circle in the sense that it can be broken into two complementary sectors containing all 
$\Alpha$'s resp. $\Beta$'s;

\item 
$\{\alpha_i\} \subset [\,0,1\,)$, 
$\{\beta_j\} \subset (\,-1,0\,]$.
\end{itemize}
To fix a scaling, set $\lambda =\exp \sum_i (\psi (\overline{\alpha_i})-\psi (\overline{\beta_i}))$ where
$\psi (x)= \frac{\Gamma '(x)}{\Gamma (x)}$
and $\overline{y}$ denotes the unique 
representative of the class $\,y \!\! \mod \Z$ in 
$(0,1]$: $\overline{y} = 1-\{-y\}$.  It follows from the multiplication formula for the gamma function that 
$\lambda \in \Q$. 
Let
$\text{univ} : U \to (\Gm \setminus \{\lambda^{-1} \})^{\text{an}}$ 
denote the universal cover. Let $\mathscr{U}$ be the weight $1-2n$  VHS whose underlying local system is constant with the fiber $\Q^n$, and the Hodge filtration is given as follows: consider the matrix $\Pi_{\mathrm{A}}(z)$ whose $j$--th column is $(z\frac{d}{dz})^j V_{\mathrm{A}}^t(S_{{\mathrm{A}}_{1}}(\lambda z), \dots, S_{\mathrm{A}_n} (\lambda z))^t $, and let $\text{Fil}^{-n/2-j} \mathcal{U}$ be the span of rows $0, \dots, j$ in $\Q^n \tensor \C$.

It is convenient to follow Deligne's and Bloch's convention and twist $\mathcal{U}$ by $\Q (1-n)$: there exists a unique weight $-1$ hypergeometric VHS $\mathcal{V}$ on 
$(\Gm \setminus \{\lambda^{-1} \})^{\text{an}}$ such that $\mathcal{U} \tensor \Q (1-n) =  \text{univ}^*  \mathcal{V} $. 
Katz's weight convention is the opposite of ours for $\mathcal{U}$.
For each $z_0 \ne \lambda^{-1}$ in $\Gm (\Q)$, his theory of $l$--adic hypergeometric sheaves 
enables one to construct naturally a weight
$(2n-1)$ {hypergeometric Galois representation} $R_{z_0}$. Finally, Magma's convention on hypergeometric motives is yet something  different: there should exist a {hypergeometric motive}  $M_{z_0}$ of weight $n-1$ such that 
$R_{z_0} = H_{\text{\'et}}( M_{z_0} \tensor \Q (-n/2),\Q_l)$ and $\mathcal{V}_{z_0} = H_\text{dR} ( M_{z_0} \tensor \Q (n/2))$. 
Conceptually, these are all minor details that affect the computations in a trivial way.   

\pph \bf Deligne's conjecture \rm (\! \cite {Roberts}, \cite{Yang2021}, and  unpublished 
computations by Candelas--de la Ossa--van Straten). With the assumptions made in the previous paragraph, it says that the value $L(M_{z_0},n/2)$ is proportional with a rational  factor to a certain minor arising from the Betti to de Rham identification for $\mathcal{V}_{z_0}$, 
or equivalently, from the period matrix for the Hodge structure
$\mathcal{V}_{z_0} \tensor \Q (-1)$. Concretely,
one expects 
$$\frac{L(M_{z_0},n/2)} {  \det \, (2 \pi \ii)^{n} \mathop{\mathrm{Re}} \Pi_{\mathrm{A}}(z_0)_{\{0, \dots, n/2-1\},\{0, \dots, n/2-1\}}} \, \in \Q,$$ where the subscript indicates the top-left quarter of the period matrix. 
Experimenting with the $L$--functions (as implemented in Magma) for the case $n=4$ corresponding to weight $3$ Calabi--Yau motives, 
one checks the identity numerically for various different values of $z_0$ for the $7$ out of the 14 MUM families that are 
non--resonant at $z=0$ (a hypergeometric differential equation is non--resonant at $0$ resp. $\infty$ if the eigenvalues $\Alpha$'s resp. $\Beta$'s of the local monodromy operator are distinct). Concretely,
the $\alpha$'s and $\beta$'s in the 7 families
are as in the left table below.

\pph \label{qua-tw} \bf The quadratic twist and the Birch--Swinnerton-Dyer period. \rm  
 Following a suggestion by Fernando Rodriguez Villegas,
we \emph{twist} the $\alpha$'s and $\beta$'s by shifting all the indices by $-\frac{1}{2}$: $\ta, \tb$  are, respectively, in the right table; $\tilde \lambda$ is now obtained from $\ta, \tb$ by the same rule as above.

\medskip

$
\begin{array}{c|c|c}
\mystrut 1 & [ \frac{1}{12},\frac{5}{12}, \frac{7}{12}, \frac{11}{12} ] & [ 0, 0, 0, 0 ]  \\ \hline \mystrut 2 & [ \frac{1}{10}, \frac{3}{10}, \frac{7}{10}, \frac{9}{10} ] & [ 0, 0, 0, 0 ]  \\ \hline \mystrut 3 & [ \frac{1}{8}, \frac{3}{8}, \frac{5}{8}, \frac{7}{8} ] & [ 0, 0, 0, 0 ]  \\ \hline \mystrut 4 & [ \frac{1}{6}, \frac{1}{4}, \frac{3}{4}, \frac{5}{6} ] & [ 0, 0, 0, 0 ]  \\ \hline \mystrut 5 & [ \frac{1}{6}, \frac{1}{3}, \frac{2}{3}, \frac{5}{6} ] & [ 0, 0, 0, 0 ]  \\ \hline \mystrut 6 & [ \frac{1}{5}, \frac{2}{5}, \frac{3}{5}, \frac{4}{5} ] & [ 0, 0, 0, 0 ]  \\ \hline \mystrut 7 & [ \frac{1}{4}, \frac{1}{3}, \frac{2}{3}, \frac{3}{4} ] & [ 0, 0, 0, 0 ] . \end{array}$
\hspace{3em}
$\begin{array}{c|c|c}
\mystrut \widetilde{1} & [ -\frac{5}{12}, -\frac{1}{12}, \frac{1}{12}, \frac{5}{12} ] & [ -\frac{1}{2}, -\frac{1}{2}, -\frac{1}{2}, -\frac{1}{2} ]  \\
\hline \mystrut \widetilde{2} & [ -\frac{2}{5}, -\frac{1}{5}, \frac{1}{5}, \frac{2}{5} ] & [ -\frac{1}{2}, -\frac{1}{2}, -\frac{1}{2}, -\frac{1}{2} ]  \\ \hline \mystrut \widetilde{3}  & [ -\frac{3}{8}, -\frac{1}{8}, \frac{1}{8}, \frac{3}{8} ] & [ -\frac{1}{2}, -\frac{1}{2}, -\frac{1}{2}, -\frac{1}{2} ]  \\ \hline \mystrut \widetilde{4}  & [ -\frac{1}{3}, -\frac{1}{4}, \frac{1}{4}, \frac{1}{3} ] & [ -\frac{1}{2}, -\frac{1}{2}, -\frac{1}{2}, -\frac{1}{2} ]  \\ \hline \mystrut \widetilde{5}  & [ -\frac{1}{3}, -\frac{1}{6}, \frac{1}{6}, \frac{1}{3} ] & [ -\frac{1}{2}, -\frac{1}{2}, -\frac{1}{2}, -\frac{1}{2} ]  \\ \hline \mystrut \widetilde{6}  & [ -\frac{3}{10}, -\frac{1}{10}, \frac{1}{10}, \frac{3}{10} ] & [ -\frac{1}{2}, -\frac{1}{2}, -\frac{1}{2}, -\frac{1}{2} ]  \\
\hline \mystrut \widetilde{7}  & [ -\frac{1}{4}, -\frac{1}{6}, \frac{1}{6}, \frac{1}{4} ] & [ -\frac{1}{2}, -\frac{1}{2}, -\frac{1}{2}, -\frac{1}{2} ] .
\end{array}$
\medskip

Put $\tilde{\mathbf{\Gamma}}_{\ta,\tb}(s) = {\tilde{\lambda}}^{1/2} \mathbf{\Gamma}_{\ta,\tb}(s)$. One has  $\tilde{\mathbf{\Gamma}}_{\ta,\tb}(s) = {\tilde{\lambda}}^{1/2} \, \mathbf{\Gamma}_{\alpha, \beta}(s+1/2).$ 
All that has been said up to now about hypergeometric Hodge structures works identically for the 7 families and 
the 7 twists. However, we expect the twist to raise the `average' analytic rank in the family. Starting with an  $L_{\ta,\tb}$ as above, 
we will construct a `biextension' variation of mixed Hodge structure formally in hypergeometric terms.  Although it is not true in general that the product of two differential operators of motivic origin is again motivic, there are situations when one can construct mixed motivic variations formally.

\pph \bf Theorem. \rm  With the assumptions
({\bf A1}) and ({\bf A2}) made
in \ref{assumptions-1}. and \ref{assumptions-2}., the differential equation 
$DL_{\ta,\tb} \, D S(z)=0$ is motivic, i.e. underlies a VMHS of geometric origin.  

\bigskip

 \bf Proof. \rm 
The idea is that under certain conditions that hold in our case we can pass from the $D$--module corresponding to a differential operator $L$ to the one corresponding to $DLD$  by successively convoluting it 
with the star resp. the shriek extension of the `constant object'~$\OO$ on $\Gm - \{1\}$
to $\Gm$. The background is \cite{Katz1990}; all references in the proof are to this book. Denote $\partial =\frac{d}{dz},  \, D= z\dd$ as above, $\, \DD = \DD_{\Gm} = \C[z,z^{-1},\dd], \,
\DD_{\A1} = \C[z,\dd].$ Let $j$ be
the open immersion $\Gm \hookrightarrow \A1$, and let inv denote the inversion map on $\Gm$. We denote the Fourier transform functor by FT.   

\proofstep
\sl Katz's lemma on indicial polynomials. \rm [2.9.5]
Write $L$ as a polynomial in $z$ whose coefficients are in turn polynomials in $D$: \:
$L=\sum_{k=0}^d z^k P_k(D)$. Then:
$P_0(y)$ has no zeroes in $\Z_{<0}$ iff
$$ \DD_{\A1} / \DD_{\A1} L  \isomorphic j_! ( \DD /  \DD L);$$
$P_0(y)$ has no zeroes in $\Z_{\ge 0}$ iff
$$  \DD_{\A1} / \DD_{\A1} L  \isomorphic j_* ( \DD /  \DD L).$$

\proofstep
\sl The D--modules
$F_k=\DDD/\DDD (D-k)$ with $k \in \Z$. \rm 
The lemma says that for $k \ge 0$, the D--module $F_k$ is isomorphic to $j_! \OO$;
for $k < 0$, the D--module $F_k$ is isomorphic to $j_* \OO$.

%

We will need a version of this: put
$E_k=\DD /\DD  (D-z(D-k))$ with $k \in \Z$. Denote by $j'$ 
the open immersion $\Gm-\{1\} \hookrightarrow \Gm$. We claim that for $k \ge 0$, the D--module $E_k$ 
is $j'_! \OO_{\Gm-\{1\} }$. For $k < 0$, the D--module $E_k$ 
is $j'_* \OO_{\Gm-\{1\} }$. Indeed, put $z=1+u$, then $D-z(D-k) = (1+u)\dd -(1+u)(1+u)\dd + (1+u)k = -(1+u) (u\dd - k)$.

\proofstep
\sl Katz's `key lemma'. \rm [5.2.3] Let the convolution sign stand for convolution with no supports on $\Gm$. 
For any holonomic module $M$ on $\Gm$ we have
\begin{equation}
	j^*\FT(j_* \mathrm{inv}_*(M)) \isomorphic M \convo \bigl( \DD / \DD (D-z)\bigr) 
\end{equation}
and [5.2.3.1]
\begin{equation}
	\mathrm{inv}_*j^*\FT(j_*M) \isomorphic M \convo \bigl( \DD / \DD (1+zD)\bigr) .
\end{equation}

\proofstep
We define the star (resp. the shriek) Ur--object to be 
$$ \bigl( \DD / \DD (1-zD)\bigr) \convo \bigl( \DD / \DD (D-z)\bigr)$$ 
resp.
$$ \bigl( \DD / \DD (1-zD)\bigr) \convo_! \bigl( \DD / \DD (D-z)\bigr).$$  

\medskip
\emph{Claim.} The star Ur--object is $E_{0}$. Proof (cf. [6.3.5]): use the key lemma
with $M = \DD / \DD (1-zD) $. The LHS becomes

	\begin{multline*}
		j^* \FT \bigl( \DD_{\A1} / \DD_{\A1} ((D+1)+z)\bigr) 
		\isomorphic
		j^* \bigl( \DD_{\A1} / \DD_{\A1} \FT((D+1)+z)\bigr)
		\\ \isomorphic  
		j^* \bigl( \DD_{\A1} / \DD_{\A1} (-D+\partial)\bigr) 
		\isomorphic 
		\DD / \DD (D-zD).
	\end{multline*}

\proofstep Let $H = P_0(D)-zP_1(D)$ be an irreducible hypergeometric operator, so that the sets of roots $\!\!\!\mod \Z$
of $P_0$ and $P_1$ are disjoint. Assume further that $P_1$ has no integer roots and $P_0$ has no integer roots in $\Z_{\ge 0}$. 
We claim that 
$$\DD / \DD H \convo E_0 \isomorphic \DD / \DD \bigl(D P_0(D-1)-  z D P_1(D-1)\bigr).$$
Indeed, in order to convolute with $E_0$ one convolutes first with 
$ \DD / \DD (D-z)$ then with  \linebreak 
$ \DD / \DD (1-zD) \isomorphic [z \mapsto -z]_* \bigl( \DD / \DD (1+zD)\bigr)$.
The result of the first convolution is simply 
$\DD / \DD ((D+1) P_0(D)-zP_1(D))$ as $P_1$ has no integer roots, [5.3.1].  
In  order to convolute with  
$\DD / \DD (1+zD)$ one now uses the second statement of the key lemma,
obtaining 
	\begin{multline*}
		\mathrm{inv}^*\bigl( \DD / \DD (-D P_0(-D-1)-\partial P_1(-D-1))\bigr) 
		\\ \isomorphic 
		\mathrm{inv}^* \bigl( \DD / \DD (-z D P_0(-D-1)-D P_1(-D-1))\bigr) 
		\\  \isomorphic
		\DD / \DD (D P_0(D-1)+ z D P_1(D-1)). 
	\end{multline*}
Finally, the effect of $[z \mapsto -z]_*$ is in simply changing the sign of 
$z$.

%

\proofstep
Let $^\vee$ denote the `passing to adjoints' anti-automorphism sending $t$ to $t$ and $\dd$ to $-\dd$, so that the formal adjoint of  
$\bigl(P_0(D-1)- z P_1(D-1)\bigr)D$ 
is $(-D-1)\bigl(P_0(-D-2)-P_1(-D-2)z\bigr).$
Assume now that $P_0$ has no integer roots.
The previous consideration applies
so convoluting  with $E_0$ we get the 
$\DD$--module corresponding to the operator
\begin{multline*}
\Bigl((-(D-1)-1)\bigl(P_0(-(D-1)-2)-P_1(-(D-1)-2)z\bigr)\Bigr)D  \\
= \,-D \bigl(P_0(-D-1)-P_1(-D-1)z\bigr)D.
\end{multline*} 
Passing to adjoints again,
$$\Bigl[-D \bigl(P_0(-D-1)-P_1(-D-1)z\bigr)D \Bigr]^\vee
=\,
-(-D-1) \bigl(P_0(D)- z P_1(D)\bigr)(-D-1)
$$
we arrive at the $\DD$--module
$$\DD / \DD \Bigl((D+1)H(D+1)\Bigr) \isomorphic \DD / \DD \Bigl(D \bigl(P_0(D-1)-  z  P_1(D-1)\bigr) D\Bigr).$$

\proofstep To finish the proof, take
$H$ to be the hypergeometric operator
whose indices are $\ta$'s and $\tb$'s shifted by $-1$, and the position of the singularities are the same. By Katz, $H$ is motivic. By the argument above, one can pass from the D-module $\DD/\DD H$ to the 
D-module  $\DD/\DD (DL_{\ta,\tb}D)$ by successively applying the motivic operations of convolution with the motivic object $E_{-1}$ and passage to duals. Hence, 
$\DD/\DD (DL_{\ta,\tb}D)$
is itself motivic, namely 
$\DD/\DD (DL_{\ta,\tb}D) \simeq \bigl((\DD/\DD L_{\ta,\tb}) \convo j_! \mathcal{O} \bigr) \convo_{\! !} j_* \mathcal{O}$.   

\qed

\pph We remark that all these considerations translate immediately  into the $l$--adic setting. We stick with Hodge modules, but what we need here is a concrete description suitable for computation. 
The significance of the twist is that the variation of mixed Hodge structure in question is a biextension VHS \cite{Hain1990}, i.e. sits in a $\Q(1) \hookrightarrow \mathcal{V} \twoheadrightarrow \Q$; this would not be the case without the twist. Think of the fiber $\mathcal{V}$ at $z_0 \in \Q$ as realized in  $H^3 ( X_{z_0}, \Q (2))$ for a threefold $X_{z_0}$. By specializing this VMHS we construct a non--trivial (in general) biextension of $H^3 ( X_{z_0}, \Q (2))$, and by relaxing the structure to a once--extension, a class in absolute Hodge cohomology $H^4_{\mathrm{Hodge}}(X_{z_0},\R (2))$. According to the Beilinson rank conjecture, this class signals the presence of a non--trivial class $c_{z_0}$ in $\mathrm{CH}^{(2)}_0(X_{z_0})\tensor \Q$.

In the language of period matrices,
in addition to the 
4 pure periods 
$$(\Phi_1(z),\Phi_2(z),\Phi_3(z),\Phi_4(z)) = (S_{\mathrm{\tilde A}_{1}}({\tilde{\lambda}} z), \dots, S_{\mathrm{\tilde A}_{4}} ({\tilde{\lambda}} z))  V_{\mathrm{\tilde A}} $$
one introduces an extension solution  
$ S_{1}( z) = \sum_{n=0}^\infty \tilde{\mathbf{\Gamma}}_{\ta,\tb}(n) z^n$ so that $DL_{\ta,\tb}S_1({\tilde{\lambda}} z)=0$ 
and the (transposed) biextension period matrix 
$$\Pi_{\tilde{\mathrm{A}}}^\text{biext} (z) = 
\bigl((z\frac{d}{dz})^{-1}, 1, \dots, (z\frac{d}{dz})^{4}\bigr)^t\; (
S_1(\tilde{\lambda} z), 
  \Phi_{1}(z)
, \Phi_{2}(z)
, \Phi_{3}(z)
, \Phi_{4}(z),
0)
$$
with the choice of the constant terms in the $0$th row being 
$$ 
\bigl(\, (1/\ta_1 +1/\ta_2)\, \tilde{\mathbf{\Gamma}}_{\ta,\tb}(0),0,0,0,0,
 (2 \pi \ii)\,
\tilde{\mathbf{\Gamma}}_{\ta,\tb}(0) \, \bigr).
$$

A version of the Birch--Swinnerton-Dyer--type conjecture \cite{Bloch1980, KontsevichZagier2001, Scholl1991}  translates into the following statement. By analogy with the elliptic curve cases two dimensions lower, one expects that  the archimedean component of the height of $c_{z_0}$  is essentially  the ratio of two minors of  
$ \mathop{\mathrm{Re}} \Pi_{\tilde{\mathrm{A}}}^\text{biext} (z_0)$:
$$h_\text{arch} (c_{z_0}) = {\tilde{\mathbf{\Gamma}}_{\ta,\tb}(0)}^{-1}\cdot \frac{\det \mathop{\mathrm{Re}} \Pi_{\tilde{\mathrm{A}}}^\text{biext} (z_0)_{\{0,1,2 \},\{0,1,2 \}}} 
{\det \mathop{\mathrm{Re}} \Pi_{\tilde{\mathrm{A}}}^\text{biext} (z_0)_{\{1,2 \},\{1,2 \}}}.$$
Assume, in addition, that the modulus $z_0 \in \Q$ is chosen so that there are no non--archimedean components of the height. Since the minor $\det {\mathop{\mathrm{Re}} \Pi_{\tilde{\mathrm{A}}}^\text{biext} (z_0)_{\{1,2 \},\{1,2 \}}}$ occurring in the denominator is nothing else but the ${\tilde{\mathbf{\Gamma}}_{\ta,\tb}(0)}^{-2}$--scaled Deligne period of $M_{z_0}$, a version of B--SD for an analytic rank 1 motive $M_{z_0}$ in a hypergeometric family as above would predict that
$$ 
r(z_0) \, := \, \frac{L'(M_{z_0},2)}{{\tilde{\mathbf{\Gamma}}_{\ta,\tb}(0)}^{-3}  \det \mathop{\mathrm{Re}} \Pi_{\tilde{\mathrm{A}}}^\text{biext} (z_0)_{\{0,1,2 \},\{0,1,2 \}} } \, \in \Q^*.
$$

\pph \bf Examples. \rm Consider the hypergeometric 
family $\widetilde{2}$ in the second table in \ref{qua-tw}. (so that ${\tilde{\mathbf{\Gamma}}_{\ta,\tb}(0)} = 32 \,(2\pi \ii)^{-4}$). One finds numerically
$$
r(1/2) \;\;  {\overset{?}{=}} \;\; -2^{3}\cdot  5^{-2} \;\;
\text{ and  } \;\;
r(1) \;\; {\overset{?}{=}} \;\;  -2^{3}\cdot 5^{-2}.
$$

More:
\bigskip

\hspace{2em} \renewcommand{\arraystretch}{1.44}
\begin{tabular}{|c|c|r|c|}
	\hline
	$\alpha$'s & $t$ & $r(t)$ & conj. value of $- t^{-3}\,r(t)$  \\ \hline 	
[-2/5, -1/5, 1/5, 2/5] & 1/7 & -0.134344023 & 1152/25 \\ \hline
[-2/5, -1/5, 1/5, 2/5] & 1/3 & -0.189629629 & 128/25 \\ \hline
[-2/5, -1/5, 1/5, 2/5] & 1/8 & -0.160000000 & 2048/25 \\ \hline
[-3/8, -1/8, 1/8, 3/8] & 1/2 & -0.0555555555 & 4/9 \\ \hline
[-3/8, -1/8, 1/8, 3/8] & 1/8 & -0.01388888888 & 64/9 \\ \hline
[-1/3, -1/4, 1/4, 1/3] & 1/6 & -\textsc{3.160493826} & 2048/3 \\ \hline
[-1/3, -1/4, 1/4, 1/3] & 1/2 & -2.666666666 & 64/3 \\ \hline
[-1/3, -1/4, 1/4, 1/3] & 1/3 & -1.580246913 & 128/3 \\ \hline
[-1/3, -1/6, 1/6, 1/3] & 1/8 & -0.1111111111 & 512/9 \\ \hline
[-1/3, -1/6, 1/6, 1/3] & 1/7 & -0.04146420466 & 128/9 \\ \hline
[-1/3, -1/6, 1/6, 1/3] & 1/2 & -0.444444444 & 32/9 \\ \hline
[-3/10, -1/10, 1/10, 3/10] & 1/2 & -0.01777777777 & 32/225 \\ \hline
[-1/4, -1/6, 1/6, 1/4] & 1/2 & -1.333333333 & 32/3 \\ \hline
\end{tabular}

\vspace{2.5em}
 \begin{center}
 \rule{5em}{1.0pt} 
 \end{center}

Much of this can be generalized to higher--rank hypergeometrics. The method can be extended to cases involving certain higher regulators as will be shown in a forthcoming paper with Matt Kerr.  
 
\bigskip

This communication is to appear in `Experimental Results'. I thank Kilian B\"onisch for checking the computations. I thank the members of the International Groupe de Travail on differential equations in Paris for many helpful discussions, and Neil Dummigan and Emre Sert\"oz for comments and corrections.  
I~thank the Max Planck Institute for Mathematics
for its hospitality during my stay there in~2021. 

I am deeply thankful to the Institut des Hautes \'Etudes Scientifiques for the extraordinary support it gave me in 2022.

\bigskip
\nocite{Beilinson1987}
\nocite{GilletSoule1984}
\nocite{Katz1990}
 
 \nocite{DoranMorgan2006}

\nocite{Bloch1980}
\nocite{BlochKato1990}
\nocite{Deligne1979}
\nocite{SwinnertonDyerBirch1975}
\nocite{GolyshevMellit2014}
\nocite{BeukersCohenMellit2015}
\nocite{RobertsVillegas2021}
\nocite{BlochdeJongSertoz2022}

\bigskip
\bigskip

\bibliographystyle{plain}




\vskip1.3cm \noindent
{\smallskip\noindent
	\sc
	ICTP Math Section \newline\noindent
	Strada Costiera 11,  Trieste 34151 Italy}

\smallskip\noindent
and 

{\smallskip\noindent \sc Algebra and Number Theory Lab \newline\noindent
Institute for Information Transmission Problems \newline\noindent
Bolshoi Karetny 19, Moscow 127994, Russia}

\smallskip\noindent \rm
golyshev@mccme.ru

\bigskip

\end{document}